\documentclass{amsart}
\usepackage{latexsym}
\usepackage{amsmath, amssymb}
\usepackage{color}
\usepackage[all]{xy}
\newtheorem{dfs}{Definition}[section]
\newtheorem{lms}[dfs]{Lemma}
\newtheorem{thms}[dfs]{Theorem}
\newtheorem{props}[dfs]{Proposition}

\newtheorem{cors}[dfs]{Corollary}
\newtheorem{rems}[dfs]{Remark}

\newtheorem*{thm*}{Theorem}

\begin{document}

\title{A universal property for the Jiang-Su algebra}
\author{Marius Dadarlat}
\address{Department of Mathematics, Purdue University,
150 N. University St., West Lafayette, IN, 47907-2067, USA}
\email{mdd@math.purdue.edu}
\author{Andrew S. Toms}
\address{Department of Mathematics and Statistics, York University,
4700 Keele St., Toronto, Ontario, Canada, M3J 1P3}
\email{atoms@mathstat.yorku.ca}
\keywords{Jiang-Su algebra, strongly self-absorbing C$^*$-algebras}
\subjclass[2000]{Primary 46L35, Secondary 46L80}

\date{\today}

\thanks{M.D. was partially supported by NSF grant \#DMS-0500693; A.T. was partially supported by NSERC}

\begin{abstract}
We prove that the infinite tensor power of a unital separable C$^*$-algebra absorbs the Jiang-Su
algebra $\mathcal{Z}$ tensorially if and only if it contains, unitally, a subhomogeneous algebra
without characters.  This yields a succinct universal property for $\mathcal{Z}$ in a category
so large that there are no unital separable C$^*$-algebras without characters known to lie
outside it.  This category moreover contains the vast majority of our
stock-in-trade separable amenable C$^{*}$-algebras,
and is closed under passage to separable superalgebras and quotients, and hence to
unital tensor products, unital direct limits, and crossed products by countable discrete groups.

One consequence of our main result is that strongly self-absorbing ASH algebras are
$\mathcal{Z}$-stable, and therefore satisfy the hypotheses of a recent
classification theorem of W. Winter.  One concludes that $\mathcal{Z}$ is the only
projectionless strongly self-absorbing ASH algebra, completing the classification of strongly
self-absorbing ASH algebras.
\end{abstract}

\maketitle

\section{Introduction}
The Jiang-Su algebra $\mathcal{Z}$ is one of the most important simple separable amenable C$^*$-algebras.
It has become apparent in recent years that the property of absorbing the Jiang-Su algebra tensorially---being
$\mathcal{Z}$-stable---is the essential property which allows one
to analyse the fine structure of a simple separable amenable C$^{*}$-algebra
using Banach algebra $\mathrm{K}$-theory and traces.  Instances of this phenomenon
can be found in the classification theory of amenable C$^{*}$-algebras (see \cite{ro1} for an introduction
to this subject, and \cite{et} for a recent survey),
and more recently in parameterisations of unitary orbits of self-adjoint operators on
Hilbert space relative to an ambient amenable C$^{*}$-algebra (see \cite{ce} and \cite{bt}).
In spite of this, $\mathcal{Z}$ has only had an {\it ad hoc} description since its discovery in
1997.  The best result concerning the uniqueness of $\mathcal{Z}$ has been that it is determined up to
isomorphism by its $\mathrm{K}$-theory and tracial simplex in a very restricted class of inductive
limit C$^{*}$-algebras.  The main result of this article concerns the $\mathcal{Z}$-stability of
the infinite tensor power of a unital separable C$^*$-algebra.   This result has several consequences for the theory of strongly-self
absorbing C$^{*}$-algebras. In particular it shows that the Jiang-Su algebra satisfies an
attractive universal property.

Let $\mathcal{C}$ be a class of unital $C^{*}$-algebras.  For $A \in
\mathcal{C}$, we let $A^{\otimes \infty}$ denote the infinite
minimal tensor product of countably many copies of $A$.  (All
C$^{*}$-algebra tensor products in the sequel are assumed to be
minimal.)  The following pair of conditions on a C$^*$-algebra $B \in
\mathcal{C}$ constitute a universal property, which we will denote by (UP) (cf.~\cite{t1}):
\begin{enumerate}
    \item[(i)] $B^{\otimes \infty} \cong B$;
    \item[(ii)] $B \otimes A^{\otimes \infty} \cong A^{\otimes
    \infty}$, $\forall A \in \mathcal{C}$.
\end{enumerate}
It is easy to see that if $B_{1}, B_{2} \in \mathcal{C}$ satisfy
properties (i) and (ii), then they are isomorphic.  To wit,
\[
B_{1} \stackrel{(i)}{\cong} B_{1}^{\otimes \infty}
\stackrel{(ii)}{\cong} B_{1}^{\otimes \infty}  \otimes B_{2}
\stackrel{(i)}{\cong} B_{1} \otimes B_{2}^{\otimes \infty}
\stackrel{(ii)}{\cong} B_{2}^{\otimes \infty}
\stackrel{(i)}{\cong} B_{2}.
\]
The question, of course, is whether a given class $\mathcal{C}$
contains an algebra $B$ satisfying (i) and (ii) at all.  Recall that
a C$^{*}$-algebra is {\it subhomogeneous} if all of its irreducible
representations are finite-dimensional, and also that
$\mathcal{Z}^{\otimes \infty} \cong \mathcal{Z}$ (\cite{key-4}).
Our main result is the following theorem.

\begin{thms}\label{main}
Let $A$ be a unital separable C$^*$-algebra.  Suppose that $A^{\otimes \infty}$ contains, unitally,
a subhomogeneous algebra without characters.  It follows that
\[
A^{\otimes\infty} \otimes \mathcal{Z} \cong A^{\otimes \infty}.
\]
In particular, $\mathcal{Z}$ satisfies (UP) in the class
$\mathcal{C}$ consisting of all such $A$.
\end{thms}
\noindent
The necessity of taking the infinite (as opposed to some finite)
tensor power of $A$ in Theorem \ref{main} is evident from the fact that $A$
itself may be subhomogeneous without characters:  any finite tensor
power of such an algebra is again subhomogenous, but no
C$^{*}$-algebra which absorbs $\mathcal{Z}$ tensorially can have
finite-dimensional representations.  This necessity also persists in
the case that $A$ is simple and infinite-dimensional:  \cite{t1}
contains an example of a simple unital separable infinite-dimensional
C$^{*}$-algebra $A$ with the property that $A^{\otimes n} \otimes
\mathcal{Z} \ncong A^{\otimes n}$ for each $n \in \mathbb{N}$ and yet
$A^{\otimes \infty} \otimes \mathcal{Z} \cong A^{\otimes \infty}$.

The bulk of the difficulty in establishing Theorem \ref{main} is
overcome by Theorem \ref{homotopic}, which states that the space of unital $*$-homomorphisms
\[
\mathrm{Hom}_{1}\left(\mathrm{I}_{p,q}; \ \mathrm{M}_{k}(\mathrm{C}(X)\right)
\equiv \mathrm{C}\left(X; \ \mathrm{Hom}_{1}(\mathrm{I}_{p,q}; \
\mathrm{M}_{k}) \right)
\]
is path connected whenever $k$ is large relative to $\mathrm{dim}(X)$;
in particular, the homotopy groups of $F^{k}:=\mathrm{Hom}_{1}(\mathrm{I}_{p,q};
\mathrm{M}_{k})$ vanish  in dimensions $<ck$ where $c\in (0,1)$ depends only on $p,q$.
We prove this result by filtering the non-manifold $F^{k}$ so that
the successive differences in the filtration are smooth manifolds,
applying Thom's transversality theorem to perturb maps into the said differences, applying a continuous
selection argument to semicontinuous fields of representations of
$\mathrm{I}_{p,q}$, and finally appealing to Kasparov's KK-theory. (Here
$\mathrm{I}_{p,q}$ denotes the prime dimension drop algebra associated
to the relatively prime integers $p$ and $q$---see Section 2 for its
definition---and $\mathrm{M}_{k}$ denotes the set of $k \times k$
matrices with complex entries.)

The subalgebra hypothesis of Theorem \ref{main} is not only
necessary ($\mathcal{Z}$ contains, unitally, a subhomogeneous algebra without
characters), but also extremely weak.  Indeed, it
is potentially vacuous for unital separable C$^*$-algebras without characters.
Examples of unital separable C$^*$-algebras
known to lie in $\mathcal{C}$ include the following (we explain why in Section 6):
\begin{enumerate}
\item[(a)] simple exact C$^{*}$-algebras containing an infinite projection;
\item[(b)] inductive limits of subhomogeneous algebras (ASH algebras) without characters;
\item[(c)] properly infinite C$^{*}$-algebras;
\item[(d)] real rank zero C$^*$-algebras without characters;
\item[(e)] C$^*$-algebras arising from minimal dynamics on a compact infinite Hausdorff space;
\item[(f)] algebras considered pathological with respect to the strong form of
Elliott's classification conjecture for separable amenable C$^*$-algebras
(\cite{ro2}, \cite{t2}, \cite{t3}).
\end{enumerate}
\noindent
The question of whether (a) and (b) together
encompass the class of simple unital separable amenable
C$^{*}$-algebras entire is an outstanding open problem.
Note that the hypotheses of Theorem \ref{main}
imply immediately that $\mathcal{C}$ is closed under passage to
unital separable superalgebras and
 quotients. In particular, $\mathcal{C}$ is closed under arbitrary unital tensor products, unital direct limits  and   crossed products by countable discrete groups.


The Jiang-Su algebra is an example of a strongly self-absorbing
C$^{*}$-algebra, i.e., a unital separable C$^{*}$-algebra
$\mathcal{D} \ncong \mathbb{C}$ such that the factor inclusion $\mathcal{D} \otimes
1_{\mathcal{D}} \hookrightarrow \mathcal{D} \otimes \mathcal{D}$ is
approximately unitarily equivalent to a $*$-isomorphism (\cite{tw}).  Such
algebras are automatically simple and amenable.  They also rare, and connected deeply
to the classification theory of separable amenable C$^*$-algebras.
Among Kirchberg algebras---simple separable amenable purely infinite C$^*$-algebras
satisfying the UCT---there are only $\mathcal{O}_2$, $\mathcal{O}_{\infty}$,
and tensor products of $\mathcal{O}_{\infty}$ with UHF algebras of ``infinite type''
(tensor products of countably many copies of a single UHF algebra);
among ASH algebras containing a non-trivial projection, we have only the
UHF algebras of infinite type.  It has been conjectured that $\mathcal{Z}$ is the only
strongly-self absorbing ASH algebra without non-trivial projections.

Recently, W. Winter proved that simple separable unital ASH algebras which are $\mathcal{Z}$-stable
and in which projections separate traces satisfy the classification conjecture of G. Elliott (\cite{w1}).
Explicitly, isomorphisms between the graded ordered $\mathrm{K}$-groups of such C$^*$-algebras can
be lifted to isomorphisms between the C$^*$-algebras.
It is known that a strongly-self absorbing C$^*$-algebra is always infinitely self-absorbing, and
that it has a unique tracial state whenever it is stably finite.
Thus, using Theorem \ref{main} and
the fact that projections trivially separate traces in a unique trace C$^*$-algebra, we
see that Winter's result
applies to strongly self-absorbing ASH algebras.
Such an algebra, when projectionless, has the same $\mathrm{K}$-theory as
$\mathcal{Z}$, and so is isomorphic to $\mathcal{Z}$ by Winter's theorem.
In other words, $\mathcal{Z}$ is, as conjectured, the only
strongly self-absorbing ASH algebra without non-trivial projections, and
this completes the classification of strongly self-absorbing ASH algebras.
We note, lest the reader find this last result too specialised, that
there are no known examples of simple unital separable amenable stably
finite C$^*$-algebras which are not ASH, and that every stably
finite strongly self-absorbing C$^*$-algebra is simple, separable, unital, and
amenable.

The sequel is organised as follows: Section 2 introduces some notation pertaining to
finite-dimensional representations of dimension drop algebras;  Section 3 examines
the homotopy groups of finite-dimensional representations of dimension drop algebras;
Section 4 provides a continuous selection theorem for sub-representations of a
semi-continuous field of representations of a dimension drop algebra over a compact
Hausdorff space;  in Section 5 we prove an extension theorem for certains maps out
of dimension drop algebras;  Section 6 combines our technical results to prove Theorem
\ref{main};  in Section 7 we make some final remarks.

\section{Preliminaries and notation.}

\subsection{Basic assumptions}
Unless otherwise noted, all morphisms in this paper are $*$-preserving algebra
homomorphisms.  We use $\mathrm{M}_k$ to denote the C$^*$-algebra of $k \times k$
matrices with complex entries.  If $X$ is a compact Hausdorff space, then $\mathrm{C}(X)$
denotes the C$^*$-algebra of continuous complex-valued functions on $X$.
If $A$ is a unital C$^*$-algebra, then $\mathcal{U}(A)$ is its unitary group.

\subsection{Finite-dimensional representations of dimension drop
algebras}\label{2.2} We assume throughout that $p$ and $q$ are
relatively prime integers strictly greater than one.  The prime
dimension drop algebra $\mathrm{I}_{p,q}$ is defined as follows:
\[
\mathrm{I}_{p,q} = \{ f \in \mathrm{C}([0,1]; \mathrm{M}_p \otimes \mathrm{M}_q) \ | \
f(0) \in \mathrm{M}_p \otimes 1_q, \ f(1) \in 1_p \otimes \mathrm{M}_q \},
\]
with the usual pointwise operations.

For any $t \in [0,1]$, we define a $*$-homomorphism $ev_t:\mathrm{I}_{p,q} \to
\mathrm{M}_{pq} \cong \mathrm{M}_p \otimes \mathrm{M}_q$ by setting $ev_t(f) = f(t)$.
If $t \in (0,1)$, then we will refer to $ev_t$ as a {\it generic evaluation}.
Now suppose that $f(0) = a \otimes 1_q$ and $f(1)=1_p \otimes b$.  Define
maps $e_0:\mathrm{I}_{p,q} \to \mathrm{M}_p$ and $e_1:\mathrm{I}_{p,q} \to \mathrm{M}_q$
by $e_0(f) = a$ and $e_1(f) = b$.  We will refer to $e_0$ and $e_1$ as {\it endpoint
evaluations}.
Set $F^k = \mathrm{Hom}_1(\mathrm{I}_{p,q}; \mathrm{M}_k)$.
 It is known that every integer $k \geq pq-p-q$
can be written as a nonnegative integral linear combination of $p$
and $q$, whence $F^{k}$ is not empty in that case.
 Each $\phi \in F^k$ is unitarily
equivalent to
\[
\tilde{\phi} = \left( \bigoplus_{i=1}^{a_{\phi}} e_0 \right) \oplus
\left( \bigoplus_{k=1}^{b_{\phi}} e_1 \right) \oplus
\left( \bigoplus_{j=1}^{c_{\phi}} ev_{x_j} \right),
\]
where $x_j \in (0,1)$.
Let us denote by $\mathrm{sp}(\phi)$ the multiset consisting of the $x_j$s.

\subsection{Spectral multiplicity}
Let $X$ be a compact Hausdorff space, and let
\[
\phi:\mathrm{I}_{p,q} \to \mathrm{C}(X) \otimes \mathrm{M}_k
\]
be a unital $*$-homomorphism.  If
$Y \subseteq X$, then $\phi|_Y$ will denote the restriction of $\phi$
to $Y$.  This restriction is actually a $*$-homomorphism into
$\mathrm{C}(Y) \otimes \mathrm{M}_k$ whenever $Y$ is closed, and is
always at least a $*$-preserving algebra homomorphism into the set of
continuous $\mathrm{M}_k$-valued functions on $Y$.

We define $N_0^{\phi}:X \to \mathbb{Z}^+$ (resp. $N_1^{\phi}:X \to \mathbb{Z}^+$) to be the
upper semicontinuous function which, at $x \in X$, returns the number of $e_0$s (resp. $e_1$s)
occurring as direct summands of $\phi|_{\{x\}}$.  Similarly, we define $N_g^{\phi}:X \to \mathbb{Z}^+$
to be the lower semicontinuous function which, at $x \in X$, returns the number of generic evaluations
occurring as direct summands of $\phi|_{\{x\}}$.  These functions are related as follows:
\begin{equation}\label{Nsum}
k = pN_0^{\phi}(x) + qN_1^{\phi}(x) + pqN_g^{\phi}(x), \ \forall x \in X.
\end{equation}

\section{Reducing the number of generic representations}\label{reduce}

Let $F_l^{k}$ denote the subset of $F^{k}$ consisting of those $\phi$ for which $\mathrm{sp}(\phi)$
contains at most $l$ points, counted with multiplicity.  It follows that the difference
$F_l^{k} \backslash F_{l-1}^{k}$ consists of those $\phi$ for which $\mathrm{sp}(\phi)$ contains
exactly $l$ points, counted with multiplicity.
Let $a$ and $b$ be nonnegative integers satisfying
\[
ap + bq +lpq = k.
\]
Let $F^k(a,b,l)$ denote the set of $\phi \in F_l^{k} \backslash
F_{l-1}^{k}$ which, up to unitary equivalence, contain
exactly $a$ direct summands of the form $e_{0}$ and $b$ direct summands
of the form $e_{1}$.  For a fixed $l$, the various $F^k(a,b,l)$ are
clopen subsets of $F_{l}^{k} \backslash F_{l-1}^{k}$.  Finally,
let $S^k(a,b,l)$ denote the subset of $F^k(a,b,l)$ consisting of
those $\phi$ which contain exactly $l$ summands of the form $ev_{1/2}$.

\begin{lms}\label{retract}
The inclusion
\[
F^k_{l-1} \stackrel{\iota}{\longrightarrow} F^k_{l-1} \cup F^k(a,b,l) \backslash S^k(a,b,l)
\]
is a deformation retract.
\end{lms}

\begin{proof}
Define a family $\{h_{t}\}_{t \in [0,1/2)}$ of continuous
self-maps of $[0,1]$ as follows:
\[
h_{t}(x) = \left\{ \begin{array}{ll} 0, & x \in [0,t] \\ (x-t)/(1-2t),
& x \in (t,1-t) \\ 1, & x \in [1-t,1] \end{array} \right. .
\]
Note that $(x,t) \mapsto h_{t}(x)$ is continuous in both variables.
Since $h_{t}$ fixes $0$ and $1$, it induces an endomorphism $\eta_t$
of $\mathrm{I}_{p,q}$ by acting on $\mathrm{Spec}(\mathrm{I}_{p,q}) \cong [0,1]$.

Let us now define a continuous map
\[
d:F^k_{l-1} \cup F^k(a,b,l) \backslash S^k(a,b,l) \to [0,1/2).
\]
On $F^k_{l-1}$, set $d = 0$.  On $F^k(a,b,l) \backslash S^k(a,b,l)$,
$d$ is the Hausdorff distance between the following two subsets of $[0,1]$:
first, the (nonempty) set of points in $\mathrm{sp}(\phi)$ which also lie
in $(0,1/2) \cup (1/2,1)$;  second, the set $\{0,1\}$.

Define a homotopy $H(t)$ of self-maps of $F^k_{l-1} \cup F^k(a,b,l) \backslash S^k(a,b,l)$
by
\[
H(t)(\phi) = \phi \circ \eta_{t d(\phi)}.
\]
Since $d(\phi) \equiv 0$ on $F^k_{l-1}$ and $\eta_0 = \mathbf{id}_{\mathrm{I}_{p,q}}$, we
have
\[
H(t)|_{F^k_{l-1}} = \mathbf{id}_{F^k_{l-1}}
\]
and
\[
H(0) = \mathbf{id}_{F^k_{l-1} \cup F^k(a,b,l) \backslash S^k(a,b,l)}.
\]
Since $\phi \circ \eta_{d(\phi)} \in F^k_{l-1}$ whenever
$\phi \in F^k(a,b,l) \backslash S^k(a,b,l)$, $\iota$ is a deformation retract.
\end{proof}

\begin{props}\label{manifold}
The topological space $F^k(a,b,l)$ can be endowed with the structure of a smooth manifold;
the subspace $S^k(a,b,l)$ is then a compact submanifold of codimension  $l^2$.
\end{props}

\begin{proof}
Note that
$M=\{\varphi:\mathrm{Hom}(C[0,1],M_{l}):Sp(\varphi)\subset(0,1)\}$
 can be naturally identified with the set of
all selfadjoint matrices in $M_{l}$ having all their eigenvalues in  $(0,1)$
and hence $M$ is homemomorphic to $\mathbb{R}^{l^{2}}$.

We are going to exhibit a free,  proper and smooth
right action of the compact Lie group
$G=U(a)\times U(b)\times U(l)$ on the manifold $X=M\times U(k)\cong \mathbb{R}^{l^{2}}\times U(k)$,
such that the quotient space is  homeomorphic to $F^{k}(a,b,l)$.
Then we invoke a result from \cite{key-1} to conclude that $X/G$ and hence $F^{k}(a,b,l)$  admits a unique smooth structure for which the quotient
map is a submersion. The uniqueness part of the same result shows that if $Y$ is a $G$-invariant submanifold
of $X$, then $Y/G$ is a submanifold of $X/G$.

If $\varphi\in M$, then $\varphi\otimes\mathrm{id}_{pq}:C[0,1]\otimes M_{pq}\to M_{l}\otimes M_{pq}$
defines by restriction a morphism on $I_{p,q}\subset C[0,1]\otimes M_{pq}$. Let us define a continuous map
$P:M\times U(k)\to F^{k}(a,b,l)$ by \[
P(\varphi,u)(f)=u[(e_0(f)\otimes1_{a})\oplus(e_1(f)\otimes1_{b})\oplus
(\varphi\otimes\mathrm{id}_{pq})(f)]u^{*},\]
where for $f\in I_{p,q}$, $f(0)=e_0(f)\otimes 1_q\in M_p \otimes M_q$ and $f(1)=1_p\otimes e_1(f)\in M_p \otimes M_q $.

One verifies that the map $P$ is surjective and that $P(\varphi,u)=P(\psi,v)$ if and only if there is
$W=(w_{0},w_{1},w)\in G$ such that $\psi=w^{*}\varphi w$ and $v=uj(W)$
where $j:G\to U(k)$ is the injective morphism\[
j(w_{0},w_{1},w)={(1}_{p}\otimes w_{0})\oplus{(1}_{q}\otimes w_{1})\oplus(w\otimes1_{pq}),\]
induced by the embedding of C{*}-algebras \[
B={(1}_{p}\otimes M_{a})\oplus(1_{q}\otimes M_{b})\oplus{(M}_{l}\otimes1_{pq})\subset M_{k}.\]
 This shows that if we define a right action of $G$ on $X=M\times U(k)$
by\[
(\varphi,u)W=(w^{*}\varphi w,vj(W)),\]
 where $W=(w_{0},w_{1},w)\in U(a)\times U(b)\times U(l)=G$,  then $P$ induces a continuous bijection
$X/G\to F^{k}(a,b,l)$. This induced map is actually a homeomorphism
since one can verify that $P$ is an open map as follows.
 Fix $(\varphi,u)$ in $X$ and
let $V$ be a neighborhood  of $(\varphi,u)$. We need to show that
 if $P(\psi,v)$
is sufficiently close to $P(\varphi,u)$, then there is $(\psi_1,v_1)$ in $V$ such
that $P(\psi_1,v_1)=P(\psi,v)$. Fix a metric $d$ for the point-norm topology
of $F^k$. Suppose that $d(P(\varphi,u),P(\psi,v))<\delta$ for some  $\delta>0$
to be specified later. Then $v^*u$ approximately commutes with the unit ball of the subalgebra
$A=(M_p\otimes 1_a) \oplus (M_q\otimes 1_b)\oplus (M_l\otimes 1_{pq})$
of $M_k$. By a classical perturbation result for finite dimensional C*-algebras,
there is a unitary $z$ in the relative commutant of $A$ in $M_k$, $z\in U(A'\cap M_k)=U(B)$ such that $\|v^*u-z\|<g(\delta)$, where $g$ is a universal positive
map with converges to $0$ when $\delta\to 0$. We can write $z=j(W)$ where
$W=(w_{0},w_{1},w)\in U(a)\times U(b)\times U(l)$, as above.
Let us set $\psi_1=w^*\psi w$ and $v_1=vj(W)$. Then $P(\psi_1,v_1)=P(\psi,v)$
and if $\delta$ is chosen sufficiently small, then $(\psi_1,v_1)$ is in $V$
since $\|u-vj(W)\|<g(\delta)$ and since $d(w^*\psi w,\varphi)\to 0$  as $\delta \to 0$, because $d(P(\varphi,1),P(\psi,u^*v))<\delta$.

Having established  the homeomorphim $F^k(a,b,l)\cong X/G$, we need to argue that
 $X/G$ is a smooth manifold.
To this purpose we apply  Proposition 5.2
on page 38 of \cite{key-1}, according to which $X/G$ is a manifold provided that $X$
is a manifold and the  action of the Lie group $G$ on $X$ is  proper,
free and smooth. Recall that a free right action $G\times X\to X$
is proper if
\begin{enumerate}
\item The set $C=\{(x,xg):x\in X,g\in G\}$ is closed in $X\times X$ and
\item The map $\iota:C\to G,$ $\iota(x,xg)=g$ is continuous.
\end{enumerate}
The first condition is easily verified if $G$ is compact as shown
in the last part of the proof of Proposition 3.1 on page 23 of \cite{key-1}.
To verify the second condition suppose that $(x_{n},x_{n}g_{n})$
converges to $(x,xg)$ in $X\times X$. Then $x_{n}$ converges to
$x$ and hence $\mathrm{dist}(xg_{n},x_{n}g_{n})=d(x,x_n)$ converges to zero.
It follows that $xg_{n}$ converges to $xg$. If $u$ is the component
of $x$ in $U(k)$, then $uj(g_{n})$ converges to $uj(g)$ in $U(k)$.
Therefore $g_{n}$ must converge to $g$.

In conclusion $F^{k}(a,b,l)$ is a manifold of dimension equal
to $\mathrm{dim}(X)-\mathrm{dim}(G)=l^{2}+k^{2}-(a^{2}+b^{2}+l^2)$.

On the other hand, $S^k(a,b,l)$ is the image in the quotient space $X/G$ of the
$G$-invariant submanifold $\{\mu\}\times U(k)$ of $X=M\times U(k)$, where $\mu(f)=f(1/2)\otimes 1_l$.
By applying \cite[Prop.~5.2]{key-1} once more we deduce that
 $S^k(a,b,l)$ is a submanifold of $F^k(a,b,l)$ of dimension  $k^{2}-(a^{2}+b^{2}+l^2)$ and hence of codimension $l^2$.
\end{proof}

\begin{props}\label{equivalence}
The inclusion $F^k_{l-1} \hookrightarrow F^k_l$ is an $(l^2-1)$-equivalence.
\end{props}

\begin{proof} All manifolds in this proof are assumed to be metrisable
and separable. Let  $S^k_l$ be the union of all $S^k(a,b,l)$.
 Let $X$ be a smooth manifold of dimension $\leq l^2-1$
and let $Y$ be a compact subset of $X$. Let $f: X \to F^k_l$
be a continuous map such that $f(Y)\cap S^k_l=\emptyset$ and let $\varepsilon >0$. We are going to show that
there is a continuous map $g: X \to F^k_l$ such that
\begin{enumerate}
 \item[(i)] $d(f(x),g(x))<\varepsilon$, for all $x\in X$,
\item[(ii)] $g=f$ on $Y$,
\item[(iii)] $g(X)\cap S^k_l=\emptyset$.
\end{enumerate}
We shall apply this perturbation result in the realm of cellular maps
from a pair of CW-complexes,
$(X,Y)=(S^n,*)$ or $(X,Y)=(S^n\times [0,1],\{*\}\times [0,1])$ to the the CW-complex $F^k_l$.
This implies that
$f$ is homotopic to $g$ via a homotopy that is constant on $Y$, assuming that
$\varepsilon$ was chosen sufficiently small.
Consequently, the inclusion
\[F^k_l\setminus  S^k_l\hookrightarrow F^k_l\]
is an $l^2-1$ equivalence. One concludes the proof by combining this fact with Lemma 3.1 

Let us turn to the proof of the pertubation result.
Note that $S^k_l$ is a compact submanifold of $F^k_l\setminus F^k_{l-1}$ and the latter
is an open subset of $F^k_l$.
We may assume that $f^{-1}(S_l^k)\neq \emptyset$, for otherwise there is nothing to prove.
By  a classical approximation result, we may assume that $f$ is  smooth
on the pre-image $f^{-1}(U)$ of some open neighborhood $U$ of $S^k_l$ in $F^k_l\setminus F^k_{l-1}$. After shinking $U$, if necessary, we can arrange that
$U\cap f(Y)=\emptyset$ since $f(Y)\cap S^k_l=\emptyset$.
Choose  a smaller open neighborhood $V$ of $S^k_l$  such that
 $\overline{V}\subset U$.

Let $f_1:f^{-1}(U)\to U$ be the restriction of $f$ to the manifold $f^{-1}(U)$.
Note that $f_1(f^{-1}(U\setminus V))$ is disjoint from $S^k_l$ by construction.
By Thom's transversality theorem \cite[p.~557]{key-3}, since $\mathrm{dim}(f^{-1}(U))\leq l^2-1$ and $S^k_l$ has codimension $l^2$ in $U$,
there is a smooth map $g_1:f^{-1}(U)\to U$ such that the image of $g_1$
is disjoint from $S^k_l$ and $g_1=f_1$ on the closed subset $f^{-1}(U\setminus V)$
of $f^{-1}(U)$ and $d(f_1(x),g_1(x))<\varepsilon$ for all $x\in f^{-1}(U)$.
Then, the map
$g$ obtained by gluing $g_1$ with the restriction of $f$ to $f^{-1}(F^k_l\setminus \overline{V})$ along $f^{-1}(U\setminus \overline{V})$ satisfies
the conditions (i)-(iii) from above.
\end{proof}

\section{Selections up to homotopy}\label{specmult}

Let $\psi:\mathrm{I}_{p,q} \to \mathrm{C}(X) \otimes \mathrm{M}_k$ be
a unital $*$-homomorphism.  Let us abuse notation slightly and write
$\psi_x$ for $\psi|_{\{x\}}$.  Set
\[
V_i = \{x \in X \ | \ N_1^{\psi}(x) = i \}
\]
and
\[
O_i = \{ x \in X \ | \ N_1^{\psi}(x) < i+1 \} = \bigcup_{j=1}^i V_j.
\]
(We will write $V_{i,\psi}$ and $O_{i,\psi}$ for $V_i$ and $O_i$, respectively,
whenever it is not clear that $\psi$ is the $*$-homomorphism with respect to which
$V_i$ and $O_i$ have been defined.)  Note that $O_i$ is open for each $i$.

Each $\psi_x$ has the form
\[
u \left( \gamma^{'} \oplus \bigoplus_{i=1}^{N_1^{\psi}(x)} e_1 \right) u^*
\]
for some $u \in \mathcal{U}(\mathrm{M}_k)$.  Set
\[
\psi^1_x = u \left( \bigoplus_{i=1}^{N_1^{\psi}(x)} e_1 \right) u^*.
\]
In words, $\psi^1_x$ is the largest direct summand of $\psi_x$
which factors through $e_1$.
We may view $\psi^1_x$ as a unital $*$-homomorphism from $\mathrm{M}_{q}$
to $\psi^1_x(1) \mathrm{M}_k \psi^1_x(1)$.  If $Y \subset X$ is closed, then we
define $\psi^1_Y:\mathrm{I}_{p,q} \to \mathrm{Map}(Y; \mathrm{M}_k)$
by
\[
\psi^1_Y(a)(x) = \psi^1_x(a), \ \ \forall x \in Y.
\]

\begin{lms}\label{restricthom}
Let $\psi:\mathrm{I}_{p,q} \to \mathrm{C}(X) \otimes \mathrm{M}_k$ be a unital
$*$-homomorphism, and let $Y \subset X$ be closed.  If $N_1^{\psi}$ is constant
on $Y$, then $\psi^1_Y$ is a $*$-homomorphism from $\mathrm{I}_{p,q}$ into
$\mathrm{C}(Y) \otimes \mathrm{M}_k \cong \mathrm{C}(Y;
\mathrm{M}_{k})$ which factors through $e_1$.
\end{lms}

\begin{proof}
Since $\phi^{1}_{Y}$ preserves pointwise operations by definition, we
need only check that $\psi^1_Y(a) \in \mathrm{C}(Y) \otimes \mathrm{M}_k$
for each $a \in \mathrm{I}_{p,q}$.  Using the fact that $N_1^{\psi}$ is constant
on $Y$ and a compactness argument, we can find some $0 < \delta < 1/2$ such that if
$ev_t$ is, up to unitary equivalence, a summand of $\psi_x$ and $x \in Y$, then
$t \in (0,1-\delta)$.  Let $\tilde{a} \in \mathrm{I}_{p,q}$ be an element
which is equal to $a$ at $1 \in \mathrm{Spec}(\mathrm{I}_{p,q})$ and
which vanishes on $[0,1-\delta) \subseteq \mathrm{Spec}(\mathrm{I}_{p,q})$.
Since both $\psi^1_Y(a)$ and $\psi(\tilde{a})|_Y$ depend only on the value of $a$ at $1$, we
conclude that they are equal.  Since $\psi(\tilde{a})|_Y \in \mathrm{C}(Y)
\otimes \mathrm{M}_k$, this completes the proof.
\end{proof}

Let $\{e_{ij}\}_{i,j=1}^{q}$ be a set of matrix units
for $\mathrm{M}_{q}$, so that $\{\psi^1_x(e_{ij})\}_{i,j=1}^{q}$ is
a set of matrix units for the image of $\psi^{1}_{x}$.  Given a subprojection
$r$ of $\psi^1_x(e_{11})$, we can generate a direct summand of $\psi^1_x$ as
follows:  for each $i,j \in \{1,\ldots,q\}$, set $f_{ij} = \psi_x^1(e_{i1}) r \psi^1_x(e_{1j})$;
notice that the $f_{ij}$ form a set of matrix units, so that
\[
\psi^1_x = \psi^{\bar{r}}_x \oplus \psi^r_x
\]
with $\psi^r_x(e_{ij}) := f_{ij}$.
We will refer to the map $\psi^r_x$ as the {\it direct summand of
$\psi^1_x$ generated by $r$}.  If $Y \subset X$ is closed and $r:Y \to \mathrm{M}_k$
is a projection-valued function such that $r(x) \leq \psi^1_x(e_{11})$ for every
$x \in Y$, then we define $\psi^r_Y:\mathrm{M}_{q} \to \mathrm{Map}(Y;
\mathrm{M}_k)$ by
\[
\psi^r_Y(a)(x) = \psi^{r(x)}_x(a).
\]
If $Y$ satisfies the hypotheses of
Lemma \ref{restricthom} and $r$ is continuous, then $\psi^r_Y$ defines a unital
$*$-homomorphism from $\mathrm{M}_{q}$ into a corner of $(\mathrm{C}(Y) \otimes
\mathrm{M}_k)$.


Let $h_{t}$ be the continuous self-map of $[0,1]$ introduced in the
proof of Lemma \ref{retract}, and let $\eta_{t}$ be the induced
endomorphism of $\mathrm{I}_{p,q}$.  Observe that $(\psi \circ
\eta_t)_x$
has fewer generic representations than $\psi_x$.  Alternatively,
\[
N_i^{\psi \circ \eta_t}(x) \geq N_i^{\psi}(x), \ \forall x \in X, \ t \in [0,1/2), \ i \in \{0,1\}.
\]
Also note that $\psi \circ \eta_t$ is homotopic to $\psi$ for each $t \in [0,1/2)$.
Since $h_{t}$ fixes
$1$, we have that
\begin{equation}\label{decomp}
(\psi \circ \eta_{t})^{1}_{x} = \psi^{1}_{x,t} \oplus
\psi^{1}_{x},
\end{equation}
for some suitable $*$-homomorphism $\psi^{1}_{x,t}$ which
factors through $e_1$.
Also, for any $0 \leq t < s < 1/2$ we have that
\begin{equation}\label{contain}
\overline{O_{i,(\psi \circ \eta_{s})}}
\subseteq O_{i,(\psi \circ \eta_{t})}.
\end{equation}
Inspection of the definition of $h_{t}$ shows that $h_{t} \circ h_{s}
= h_{s^{'}}$ for some $s^{'} \in [0,1/2)$.


\begin{lms}\label{extend}
Let $\psi:\mathrm{I_{p,q}} \to \mathrm{C}(X) \otimes \mathrm{M}_k$ be a unital $*$-homomorphism,
and let $i \in \{1,\ldots,\lfloor k/q \rfloor - 1\}$ be given.  Suppose that the following statements hold:
\begin{enumerate}
\item[(i)] there is a continuous and constant rank projection-valued map
\[
Q: \overline{O_{i,\psi}} \to \mathrm{M}_k
\]
corresponding to a trivial vector bundle ($O_{i,\psi}$ is assumed to
be nonempty);
\item[(ii)] for each $x \in \overline{O_{i,\psi}}$,
\[
Q(x) \leq \psi^{1}_{x}(e_{11})
\]
and
\[
\mathrm{rank}(Q(x)) + 2 \mathrm{dim}(X) \leq \mathrm{rank}(\psi^{1}_{x}(e_{11}));
\]
\item[(iii)] the map $\psi^{Q(x)}_{\overline{O_{i,\psi}}}$ defines a $*$-homomorphism from
$\mathrm{M}_{q}$ to $\mathrm{C}(\overline{O_{i,\psi}}) \otimes \mathrm{M}_k$.
\end{enumerate}
Then, there is a unital $*$-homomorphism
$\gamma:\mathrm{I}_{p,q} \to \mathrm{C}(X) \otimes \mathrm{M}_k$ homotopic to $\psi$
such that the following statements hold:
\begin{enumerate}
\item[(i)] there is a continuous and constant rank projection-valued map
\[
\tilde{Q}:\overline{O_{i+1,\gamma}} \cup \overline{O_{i,\psi}} \to \mathrm{M}_k
\]
which corresponds to a trivial vector bundle and extends the function $Q$ above;
\item[(ii)] for each $x \in \overline{O_{i+1,\gamma}} \cup \overline{O_{i,\psi}}$,
\[
\tilde{Q}(x) \leq \gamma^{1}_{x}(e_{11})
\]
and
\[
\mathrm{rank}(\tilde{Q}(x)) + 2 \mathrm{dim}(X) \leq \mathrm{rank}(\gamma^{1}_{x}(e_{11}));
\]
\item[(iii)] the map $\tilde{\phi}:\mathrm{M}_{q} \to \mathrm{Map}(\overline{O_{i+1,\gamma}}
\cup \overline{O_{i,\psi}}; \mathrm{M}_k)$
which agrees with $\psi^{Q(x)}_{\overline{O_{i,\psi}}}$ on $\overline{O_{i,\psi}}$ and is equal to
$\gamma^{\tilde{Q}(x)}_x$ otherwise is in fact equal to
$\gamma^{\tilde{Q}(x)}_{\overline{O_{i+1,\gamma}}
\cup \overline{O_{i,\psi}}}$, and is a $*$-homomorphism into
$\mathrm{C}(\overline{O_{i+1,\gamma}} \cup \overline{O_{i,\psi}})
\otimes \mathrm{M}_{k}$.
\end{enumerate}
\end{lms}

\begin{proof}
Choose $t \in (0,1/2)$, and set $\gamma = \psi \circ \eta_{t}$.
This ensures that $\psi$ and $\gamma$ are homotopic.
Note that $A := \overline{O_{i+1,\gamma}} \backslash O_{i,\psi}$ is closed
in $X$, and that $Q$ is already defined on $O_{i,\psi}$.  The task of
extending $Q$ to $\tilde{Q}$ is thus
reduced to the problem of extending $Q|_{A \cap \overline{O_{i,\psi}}}$ to all of $A$
and proving that our extension satisfies (i), (ii), and (iii) in the conclusion
of the lemma.

Any constant rank extension of $Q|_{A \cap \overline{O_{i,\psi}}}$ to
all of $A$ will automatically satisfy the rank requirement in conclusion (ii).
So, to begin, let us extend $Q|_{A \cap \overline{O_{i,\psi}}}$ to a trivial constant
rank projection $\tilde{Q}$ defined on all of $A$ and subordinate to
$\gamma^{1}_{x}(e_{11})$ at each $x \in A$.
For each $x \in A$, we have that
\begin{equation}\label{rank}
\mathrm{rank}(\gamma^{1}_{x}(e_{11})) \geq
\mathrm{rank}(\psi^{1}_{x}(e_{11})) \geq 2 \mathrm{dim}(X) + \mathrm{rank}(Q).
\end{equation}
By (\ref{contain}), $A$ is a closed subset of
$V_{i+1} = O_{i+1,\psi} \backslash O_{i,\psi}$.  In particular,
the rank of $\psi^{1}_{x}(e_{11})$ is constant on $A$, and
$x \mapsto \psi^{1}_{x}(e_{11})$ is a continuous and constant rank
projection valued function on $A$, call it $R(x)$.  Moreover,
$Q(x) \leq R(x) \leq \gamma^1_x(e_{11})$ for each $x \in A \cap
\overline{O_{i,\psi}}$.
It is a general fact
that $Q|_{A \cap \overline{O_{i,\psi}}}$ can be extended to a trivial projection defined
on $A$ and subordinate to $R$, as desired---all that is required
for this is the rank inequality of (\ref{rank}).
Our $\tilde{Q}$ thus satisfies parts (i) and (ii) in
the conclusion of the lemma.

Let us now establish part (iii) in the conclusion of the lemma.  Observe that
at each point in $X$, $\gamma^1_x$ decomposes as the direct sum of
$\psi^1_x$ and a second morphism, say $\lambda$ (cf. (\ref{decomp}) above).
Also recall that $Q(x) = \tilde{Q}(x)$ for each $x \in A \cap \overline{O_{i,\psi}})$.
It follows that for any $x \in A \cap \overline{O_{i,\psi}}$ and $i, j \in \{1,\ldots,q\}$, we have
\begin{eqnarray*}
\gamma^{\tilde{Q}(x)}_x(e_{ij}) & = & \gamma^1_x(e_{i1}) \tilde{Q}(x) \gamma^1_x(e_{1j}) \\
& = & (\psi^1_x(e_{i1}) \oplus \lambda(e_{i1})) Q(x) (\psi^1_x(e_{1j}) \oplus \lambda(e_{1j})) \\
& = & \psi^1_x(e_{i1}) Q(x) \psi^1_x(e_{1j}) \\
& = & \psi^{Q(x)}_x(e_{ij}),
\end{eqnarray*}
and so
\[
\gamma^{\tilde{Q}(x)}_x = \gamma^{Q(x)}_{x} = \psi^{Q(x)}_x, \ \forall
x \in A \cap \overline{O_{i,\psi}}.
\]
Our problem is thus reduced to proving that $\gamma^{\tilde{Q}}_A$ is a
$*$-homomorphism from $\mathrm{M}_{q}$ into $\mathrm{C}(A) \otimes
\mathrm{M}_{k}$.  In a near repeat of our calculation above we have
\begin{eqnarray*}
\gamma^{\tilde{Q}(x)}_x(e_{ij}) & = & \gamma_x^1(e_{i1}) \tilde{Q}(x) \gamma^1_x(e_{1j}) \\
& = & (\psi^1_x(e_{i1}) \oplus \lambda(e_{i1})) \tilde{Q}(x) (\psi^1_x(e_{1j}) \oplus \lambda(e_{1j})) \\
& = & \psi^1_x(e_{i1}) \tilde{Q}(x) \psi^1_x(e_{1j}) \\
& = & \psi^{\tilde{Q}(x)}_x(e_{ij});
\end{eqnarray*}
$\psi^{\tilde{Q}(x)}_x(e_{1j})$ is defined because $\tilde{Q}(x) \leq
R(x) \equiv \psi^1_x(e_{11})$
by construction.  We conclude that $\psi^{\tilde{Q}}_A =
\gamma^{\tilde{Q}}_A$.  Since $A \subseteq V_{i+1}$,
the function $N_1^{\psi}$ is constant on $A$.  It now follows
from Lemma \ref{restricthom} that
$\psi^{\tilde{Q}}_{A}$ is a $*$-homomorphism, completing the proof.
\end{proof}

\begin{props}\label{fdsplit}
Let $X$ be a compact metric space of covering dimension $d < \infty$, and let
$\psi:\mathrm{I}_{p,q} \to \mathrm{C}(X) \otimes \mathrm{M}_{k}$ be a
unital $*$-homomorphism.  Suppose that $N_1^{\psi}(x) > 2d$ for
each $x \in X$, and let $m$ be the miminum value taken by $N_1^{\psi}$ on $X$.
Then, $\psi$ is homotopic to a unital $*$-homomorphism
$\gamma:\mathrm{I}_{p,q} \to \mathrm{C}(X) \otimes \mathrm{M}_{k}$ with
the following property: $\gamma$
can be decomposed as a direct sum $\eta \oplus \phi$, where $\phi$ is unitarily
equivalent to
\[
\bigoplus_{j=1}^{m-2d} e_1
\]
inside $\mathrm{C}(X) \otimes \mathrm{M}_{k}$.
\end{props}

\begin{proof}
We will proceed by iterating Lemma \ref{extend}.  First, we
require some initial data satisfying the hypotheses of the said lemma.
Suppose that $O_{j,\psi}$ is empty for $0 \leq j \leq i$,
and that $O_{i+1,\psi}$ is not empty.  It follows that for some
choice of $t \in (0,1/2)$, the morphism $\Delta := \psi \circ
\eta_{t}$  has $O_{i+1,\Delta}$ nonempty.
Since $O_{i,\psi}$ is empty, we have that
\[
V_{i+1,\psi} = O_{i+1,\psi} \supseteq \overline{O_{i+1,\Delta}}.
\]
In particular, $N_1^{\psi}$ is constant on
$\overline{O_{i+1,\Delta}}$, and so
$\psi^{1}_{\overline{O_{i+1,\Delta}}}$ is a $*$-homomorphism.
The map $x \mapsto \psi^{1}_{x}(e_{11})$ is continuous and
projection-valued on $\overline{O_{i+1,\Delta}}$.  Using the
stability properties of vector bundles, we may find a continuous
and projection-valued map $Q:\overline{O_{i+1,\Delta}} \to
\mathrm{M}_{k}$ which is subordinate to $\phi^{1}_{x}(e_{11})$
at each $x \in \overline{O_{i+1,\Delta}}$, has constant rank equal
to $m-2d$, and corresponds to a trivial vector bundle.
The maps $\Delta$, $Q$, and $\Delta^{Q}_{\overline{O_{i+1,\Delta}}}$
thus constitute acceptable initial data for Lemma \ref{extend}.

Notice that in part (iii) of the conclusion of Lemma \ref{extend},
we may replace the set $\overline{O_{i+1,\gamma}}
\cup \overline{O_{i,\psi}}$ with the smaller set
$\overline{O_{i+1,\gamma}}$.
Beginning with the initial data constructed above, we iterate this
modified version of Lemma \ref{extend}
as many times as is necessary in order to arrive at a
$*$-homomorphism $\gamma$, homotopic to $\psi$ by construction,
which has the property that $O_{l+1,\gamma} = X$.  This map
has the desired direct summand $\phi$ upon restriction to
$\overline{O_{l,\gamma}}$, as provided by part (iii) in the conclusion
of Lemma \ref{extend}.  In order to extend $\phi$ to all of $X$,
simply follow the proof of Lemma \ref{extend} with $A:=O_{l+1,\gamma}
\backslash O_{l,\gamma} = V_{l+1,\gamma}$.  The proof works because
this choice of $A$ is closed.

The unitary equivalence of $\phi$ with
\[
\bigoplus_{j=1}^{m-2d} e_1
\]
follows from two facts:  first, the images of the projections $e_{11},e_{22},\ldots$ under $\phi$
all correspond to trivial vector bundles by construction; second, the complement of the sum of
these images has rank larger than $d$ and the same $\mathrm{K}_0$-class as a trivial vector
bundle, whence it is unitarily equivalent to any projection corresponding to a trivial vector
bundle of the same rank.
\end{proof}

\begin{rems}\label{switch} {\rm By replacing $1$ with $0$ and $p$ with $q$ in this
subsection, Proposition \ref{fdsplit} can be restated with $N_0^{\psi}$
substituted for $N_1^{\psi}$ and $e_0$ substituted for $e_1$.  This fact
will be used in the proof of the Lemma \ref{maxgen} below. }
\end{rems}

\begin{lms}\label{maxgen}
Let $X$ be a compact metric space of covering dimension $d < \infty$.  Let there be given a unital
$*$-homomorphism $\phi:\mathrm{I}_{p,q} \to \mathrm{C}(X) \otimes \mathrm{M}_k$ with the property
that $N_g^{\phi} < L$ on $X$ for some $L > 2d$. Assume that $k > (2L+2d+1)pq$. It follows that
$\phi$ is homotopic to a second morphism $\psi$ with the property that
\[
N_1^{\psi} \geq \left\lfloor \frac{k-2pq(L+d)}{q} \right\rfloor.
\]
\end{lms}

\begin{proof}

Since $N_0^{\phi}$ is upper semicontinuous, the set
\[
F_z^{\phi} := \{x \in X \ | \ N_0^{\phi}(x) \geq z\}
\]
is closed for every $z \in \mathbb{Z}^+$.  We also have that $F_z^{\phi} \subseteq
F_{z^{'}}^{\phi}$ whenever $z \geq z^{'}$.

\vspace{3mm}
\noindent
{\bf Claim 1:}   $F_{z+qL}^{\phi}$ is contained in the interior  of
$F_z^{\phi}$.

\vspace{1mm} \noindent {\it Proof of claim.} Let $x \in F_{z+qL}$. By the upper semicontinuity of
$N_0^{\phi}$ and $N_1^{\phi}$, there is a neighborhood $V$ of $x$ such that
$N_0^{\phi}(x)-N_0^{\phi}(v)\geq 0$ and $N_1^{\phi}(x)-N_1^{\phi}(v)\geq 0$ for all $v\in V$.
 Using (\ref{Nsum}) and the assumption that $N^{\phi}_g < L$ we have that
 \begin{eqnarray*}p (N_0^{\phi}(x)-N_0^{\phi}(v))&\leq& p (N_0^{\phi}(x)-N_0^{\phi}(v))+q(N_1^{\phi}(x)-N_1^{\phi}(v))\\&=&
pq(N_g^{\phi}(v)-N_g^{\phi}(x))< pqL
 \end{eqnarray*}
for all $v \in V$.  It follows that $N_0^{\phi}(v)> N_0^{\phi}(x)-qL\geq z+qL-qL=z$ and hence
$V\subset F_z^{\phi}$.
 \hfill $\Box$

\vspace{3mm}
\noindent
>From the claim above we conclude that for each $n \in \mathbb{N}$, there is an open set
$U_{nqL}^{\phi} \subseteq X$ such that $ F_{nqL}^{\phi} \subseteq U_{nqL}^{\phi}
\subseteq F_{(n-1)qL}^{\phi}$.

Consider the following chain of inclusions:
\[
F_{qL}^{\phi} \supseteq U_{2qL}^{\phi} \supseteq F_{2qL}^{\phi} \supseteq
U_{3qL}^{\phi} \supseteq F_{3qL}^{\phi} \supseteq \cdots.
\]
It follows from the definition of $\eta_s$ that each $x \in
F^{\phi}_z$ is an interior
point of $F_z^{\phi \circ \eta_s}$ whenever $s > 0$.  We may thus assume, by modifying
$U_{nqL}^{\phi}$ is necessary, that for some $s_0 \in
[0,1/2)$ and each $n \in \mathbb{N}$ we have
\[
F_{nqL}^{\phi \circ \eta_{s_0}} \supseteq \overline{U_{nqL}^{\phi}} \supseteq U_{nqL}^{\phi}
\supseteq F_{nqL}^{\phi}.
\]
Set $\phi_0=\phi \circ \eta_{s_0}$. Let us consider the finite set $S$ consisting of those
positive integers $n$ such that $nL>2d$ and $F_{nqL}^{\phi}\neq \emptyset$. Since $L > 2d$ by
hypothesis, we have $1 \in S$ whenever $S$ is nonempty.

\vspace{3mm} \noindent {\bf Claim 2:}  There is a $t_{0} \in [0,1/2)$ such that the following
statements hold for $\phi^{{'}} := \phi_{0} \circ \eta_{t_{0}}$:
\begin{enumerate}
\item[(i)] for every $n\in S$, the restriction $\phi^{'}|_{\overline{U_{nqL}^{\phi}}}$
has a direct summand of the form
\[
\gamma_{n}:=\bigoplus_{i=1}^{(nL-2d)q} e_0 = \bigoplus_{j=1}^{nL-2d}
ev_0;
\]
\item[(ii)] $\gamma_n|_{\overline{U_{mqL}^{\phi}}}$ is a direct summand of
$\gamma_m$ whenever $m > n$.
\end{enumerate}

\vspace{1mm} \noindent {\it Proof of claim.} We proceed by induction on $n$.  Let $n=1$. Upon
adjusting Proposition \ref{fdsplit} according to Remark \ref{switch}, we may apply it to
$\phi_{0}|_{\overline{U_{qL}^{\phi}}}$ and find a unital $*$-homomorphism
\[
\phi_{1}:\mathrm{I}_{p,q} \to \mathrm{C}(\overline{U_{qL}^{\phi}}) \otimes \mathrm{M}_k
\]
homotopic to $\phi_0|_{\overline{U_{qL}^{\phi}}}$ such that $\phi_{1}$ has a direct summand of the
form
\[
\gamma_{1}:=\bigoplus_{i=1}^{(L-2d)q} e_0 = \bigoplus_{j=1}^{L-2d} ev_0.
\]
The homotopy which arises in the proof of Proposition \ref{fdsplit} comes from composing
$\phi_0|_{\overline{U_{qL}^{\phi}}}$ with a path of $\eta_t$s. This homotopy may be extended to
all of $X$ by composing $\phi_0$ with the same path of $\eta_t$s, and so we may assume that
$\phi_{1}$ is defined on all of $X$.

Now let us suppose that we have found a unital $*$-homomorphism $\phi_{n}: \mathrm{I}_{p,q} \to
\mathrm{C}(X) \otimes \mathrm{M}_{k}$ homotopic to $\phi_{0}$ via composition with $\eta_{t}$s
such that the following statements hold:
\begin{enumerate}
\item[(i)] for every $k \in \{1,\ldots,n\} \subseteq S$,
the restriction $\phi_{n}|_{\overline{U_{kqL}^{\phi}}}$ has
a direct summand of the form
\[
\gamma_{k}:=\bigoplus_{i=1}^{(kL-2d)q} e_0 = \bigoplus_{j=1}^{kL-2d}
ev_0;
\]
\item[(ii)] $\gamma_k|_{\overline{U_{mqL}^{\phi}}}$ is a direct summand of
$\gamma_m$ whenever $m > k$.
\end{enumerate}
Assuming that $n+1\in S$, we will construct $\phi_{n+1}$, homotopic to $\phi_{n}$ via composition
with $\eta_{t}$s, such that the statements (i) and (ii) above hold with $n$ replaced by $n+1$.
Through successive applications of this inductive step we will arrive at the map $\phi^{'}$
required by the claim.

If $\overline{U_{(n+1)qL}^{\phi}}=\emptyset$ then we simply put $\phi^{'} = \phi_{n}$; suppose
that $\overline{U_{(n+1)qL}^{\phi}} \neq \emptyset$.  Let $p_n$ be the image of the unit of
$\mathrm{I}_{p,q}$ under $\gamma_{n}$, restricted to $\overline{U_{(n+1)qL}^{\phi}}$.  Applying
Proposition \ref{fdsplit} to the cut-down map
\[
(1-p_n)( \phi_{n}|_{\overline{U_{(n+1)qL}^{\phi}}}) (1-p_n),
\]
we find a unital $*$-homomorphism
\[
\phi_{n+1}:\mathrm{I}_{p,q} \to
\mathrm{C}(\overline{U_{(n+1)qL}^{\phi}}) \otimes \mathrm{M}_{k}
\]
(which, as in the establishment of the base case, arises from composing
$\phi_{n}|_{\overline{U_{(n+1)qL}^{\phi}}}$ with some $\eta_{t}$)
admitting a direct summand
\[
\alpha_{n+1} = \bigoplus_{i=1}^{Lq} e_0 = \bigoplus_{j=1}^{L} ev_{0}.
\]
Define
\[
\gamma_{n+1} = \alpha_{n+1} \oplus \gamma_{n}|_{\overline{U_{(n+1)qL}^{\phi}}}.
\]
As before, we may extend the definition of $\phi_{n+1}$ to all of $X$. These choices of
$\phi_{n+1}$ and $\gamma_{n+1}$ establish the induction step of our argument, proving the claim.
\hfill $\Box$

\vspace{3mm} Set $\alpha_1 = \gamma_1$.  Choose, for each $n\in S$, a continuous map
$g_n:\overline{U_{nqL}^{\phi}} \to [0,1]$ with the property that $g_n$ is identically zero off
$U_{nqL}^{\phi}$ and identically one on $F_{nqL}^{\phi}$. Notice that at any given $x \in X$, at
most one of the $g_n$s defined at $x$ can take a value other than one.

Define a homotopy $H:[0,1] \to \mathrm{Hom}_1(\mathrm{I}_{p,q}; \mathrm{C}(X) \otimes \mathrm{M}_k)$
as follows:
\begin{enumerate}
\item[(i)] $H(0) = \phi^{'}$.
\item[(ii)] For $t \in (0,1]$, $n$ such that $\overline{U_{nqL}^{\phi}} \backslash
\overline{U_{(n+1)qL}^{\phi}} \neq \emptyset$, and $x \in \overline{U_{nqL}^{\phi}} \backslash
\overline{U_{(n+1)qL}^{\phi}}$, $H(t)_x$ differs from $\phi^{'}_x$ as follows:  with
$\gamma_n^{\perp}$ denoting the complement of $\gamma_n$ inside
$\phi^{'}|_{\overline{U_{nqL}^{\phi}}}$ we have
\begin{eqnarray*}
\phi^{'}_x & = & (\gamma_n^{\perp})_x \oplus \bigoplus_{i=1}^{n} (\alpha_i)_x \\
& = & (\gamma_n^{\perp})_x \oplus  (\oplus_{j=1}^{L-2d} ev_{0}) \oplus
\underbrace{(\oplus_{j=1}^{L} ev_{0}) \oplus \cdots \oplus (\oplus_{j=1}^{L} ev_{0})}_{\mathrm{n-1
\ times}},
\end{eqnarray*}
where the summands of the form $\oplus_{j=1}^{L} ev_{0}$ correspond, in order, to the
$(\alpha_i)_x$s with $i > 1$;  on the other hand, modifying the second equation above, we set
\[ H(t)_x =
(\gamma_n^{\perp})_x \oplus (\oplus_{j=1}^{L-2d} ev_{t g_{1}(x)})\oplus (\oplus_{j=1}^{L} ev_{t
g_{2}(x)}) \oplus \cdots \oplus (\oplus_{j=1}^{L} ev_{t g_n(x)}).
\]
\item[(iii)] for $t \in (0,1]$ and $x \notin \overline{U_{qL}^{\phi}}$, we set $H(t)_x = \phi^{'}_x$.
\end{enumerate}

\vspace{3mm} \noindent To see that $H(t)$ is a homotopy, one need only check the continuity of
$H(t)_x$ at the boundary of $\overline{U_{nqL}^{\phi}}$.  This amounts to checking that if $x,y
\in X$ are close and $x \in \partial \overline{U_{nqL}^{\phi}}$, then $H(t)_x$ and $H(t)_y$ are
close.  Each $y$ sufficiently close to $x$ is either in $\overline{U_{nqL}^{\phi}} \backslash
\overline{U_{(n+1)qL}^{\phi}}$ or $\overline{U_{(n-1)qL}^{\phi}} \backslash
\overline{U_{nqL}^{\phi}}$ by Claim 1, so we need only address this situation. If $y \in
\overline{U_{nqL}^{\phi}} \backslash \overline{U_{(n+1)qL}^{\phi}}$, then $H(t)_{x}$ and $H(t)_y$
are close by part (ii) of the definition of $H(t)$ and the fact that $(\gamma_n^{\perp})_x$ is
continuous in $x$ on $\overline{U_{nqL}^{\phi}} \backslash \overline{U_{(n+1)qL}^{\phi}}$ (it is
the complement of $(\gamma_n)_x$, and the latter is continuous in $x$ on
$\overline{U_{nqL}^{\phi}} \backslash \overline{U_{(n+1)qL}^{\phi}}$ by construction).  If, on the
other hand, $y \in \overline{U_{(n-1)qL}^{\phi}} \backslash \overline{U_{nqL}^{\phi}}$, then we
must check that
\[ H(t)_x = (\gamma_n^{\perp})_x \oplus (\oplus_{j=1}^{L-2d} ev_{t
g_1(x)})\oplus (\oplus_{j=1}^{L} ev_{t g_{2}(x)}) \oplus \cdots \oplus (\oplus_{j=1}^{L} ev_{t
g_n(x)})
\]
is close to
\[ H(t)_y = (\gamma_{n-1}^{\perp})_y \oplus (\oplus_{j=1}^{L-2d} ev_{t g_1(y)})\oplus
(\oplus_{j=1}^{L} ev_{t g_{2}(y)}) \oplus \cdots \oplus (\oplus_{j=1}^{L} ev_{t g_{n-1}(y)})
\]
Since $\oplus_{j=1}^{L} ev_{t g_{k}(x)}$ is continuous on $ \overline{U_{(n-1)qL}^{\phi}}$ for
each $k \in \{1,\ldots,n-1\}$, we need only check that $(\gamma_n^{\perp})_x \oplus
(\oplus_{j=1}^{L} ev_{t g_n(x)})$ is close to $(\gamma_{n-1}^{\perp})_y$.  By the definition of
$g_n$, we have $g_n(x)=0$, so that
\[
(\gamma_n^{\perp})_x \oplus (\oplus_{j=1}^{L} ev_{t g_n(x)}) = (\gamma_n^{\perp})_x \oplus
(\oplus_{j=1}^{L} ev_{0}) = (\gamma_{n-1}^{\perp})_x.
\]
Our desired conclusion now follows from the continuity of $\gamma_{n-1}^{\perp}$ on $
\overline{U_{(n-1)qL}^{\phi}}$.

Put $\psi = H(1)$, so that $\psi$ is homotopic to $\phi$, as required. To complete the proof of
the lemma, we analyse the function $N_1^{\psi}$. The homotopy $H$ leaves untouched those direct
summands of $\phi^{'}_x$ of the form $e_1$, whence $N_1^{\psi} \geq N_1^{\phi^{'}}$. By
construction, we have that $N_1^{\phi^{'}}(x) \geq N_1^{\phi}(x)$ for each $x \in X$.

First consider the case $S =\emptyset$.  (Note that this assumption implies $F_{qL}^{\phi} =
\emptyset$.)  It follows that $N_0^{\phi} < qL$ on $X$. Using (\ref{Nsum}) we obtain the following
for each $x \in X$:
\begin{eqnarray*}
qN_1^{\phi}(x)& = & k-pN_0^{\phi}(x)- pqN_g^{\phi}(x)\\ & >&k-pqL-pqL\geq k-2pq(L+d).
\end{eqnarray*}
We conclude that
\[
N_1^{\psi}(x) \geq N_1^{\phi^{'}}(x) \geq N_1^{\phi}(x) > (k-2pq(L+d))/q, \] as desired.

Now suppose that $S \neq \emptyset$ and write $S = \{1,\ldots,n\}$.  We have
\[
X = X\backslash \overline{U_{qL}^{\phi}} \cup \left( \bigcup_{s=1}^{n-1} \overline{U_{sqL}^{\phi}}
\backslash \overline{U_{(s+1)qL}^{\phi}} \right) \cup \overline{U_{nqL}^{\phi}}.
\]
Fix $x \in X$. In light of the partition of $X$ above, we consider two cases, depending on whether
or not $x \in \overline{U_{qL}^{\phi}}$.

First suppose that $x \in X \backslash \overline{U_{qL}^{\phi}}$.  In this case the definition of
the homotopy $H$ implies that $\psi_x =(\phi \circ \eta_t)_x$ for some $t \in [0,1)$.  By the
spectral properties of $\eta_t$ we conclude that $N_1^{\psi}(x) \geq N_1^{\phi}(x)$.  From here,
one simply applies the argument from the $S=\emptyset$ case.

Second, suppose that $x \in \overline{U_{sqL}^{\phi}} \backslash \overline{U_{(s+1)qL}^{\phi}}$
for some $s \in \{1,\ldots,n\}$, with the convention that $\overline{U_{(n+1)qL}^{\phi}} =
\emptyset$.  From the definition of the homotopy $H$ we have
\begin{eqnarray*}
\psi_x &=&\left[(\gamma_s^{\perp})_x \oplus (\oplus_{j=1}^{L-2d} ev_{t g_{1}(x)})\oplus
(\oplus_{j=1}^{L} ev_{t g_{2}(x)}) \oplus \cdots \oplus (\oplus_{j=1}^{L} ev_{t g_s(x)})
\right]_{t=1}
\\
&=&(\gamma_s^{\perp})_x \oplus \left( \bigoplus_{j=1}^{(s-1)L-2d} ev_1 \right) \oplus
(\oplus_{j=1}^{L} ev_{g_s(x)})
\end{eqnarray*}
(the last line follows from the fact that $g_k$ is by definition identically one on
\[
\overline{U_{sqL}^{\phi}} \subseteq F_{(s-1)qL}^{\phi}
\]
whenever $k \leq s-1$).  By construction we have $(\gamma_s^{\perp})_x \oplus (\gamma_s)_x = (\phi
\circ \eta_t)_x$ for some $t \in [0,1)$.  This, by the spectral properties of $\eta_t$, implies
that
\[
N_1^{(\gamma_s^{\perp})_x \oplus (\gamma_s)_x}(x) \geq N_1^{\phi}(x).
\]
All of the $e_1$ summands of $(\gamma_s^{\perp})_x \oplus (\gamma_s)_x$ are contained in
$(\gamma_s^{\perp})_x$, so we conclude that $N_1^{(\gamma_s^{\perp})_x}(x) \geq N_1^{\phi}(x)$.
Combining this with our formula for $\psi_x$ above, we have
\begin{equation}\label{bound}
N_1^{\psi}(x) \geq N_1^{\phi}(x) + p[(s-1)L-2d].
\end{equation}
Since $F_{(s+1)qL}^{\phi} \subseteq \overline{U_{(s+1)qL}^{\phi}}$, we have $x \in X \backslash
F_{(s+1)qL}^{\phi}$.  In particular, $N_0^{\phi}(x) < (s+1)qL$.  Combining this last fact with
(\ref{Nsum}) and the inequality $N_g^{\phi} < L$ yields
\[
N_1^{\phi}(x) > \frac{k-(s+2)pqL}{q}.
\]
Combining the inequality above with (\ref{bound}) then yields
\[
N_1^{\psi}(x) \geq \frac{k-(s+2)pqL}{q}+p(sL-2d) = \frac{k-2pq(L+d)}{q},
\]
as desired.

\end{proof}
\begin{cors}\label{fdsplit2}
Let $X$ be a compact metric space of covering dimension $d < \infty$.  Let there be given a unital
$*$-homomorphism $\phi:\mathrm{I}_{p,q} \to \mathrm{C}(X) \otimes \mathrm{M}_k$ with the property
that $N_g^{\phi} < 2d+1$ on $X$. Assume that $k > (8d+4)pq$.  It follows that $\phi$ is homotopic
to a second $*$-homomorphism $\psi$ which has a direct summand of the form
\[
 \bigoplus_{i=1}^{M} ev_{1/2},
\]
where
\[
M = \left\lfloor \frac{k-pq(8d+4)}{pq} \right\rfloor.
\]
\end{cors}

\begin{proof}
Apply Lemma \ref{maxgen} to $\phi$ with $L=2d+1$.  This yields a map $\phi^{'}$ homotopic to
$\phi$ with the property that
\[
N_1^{\phi^{'}} \geq \left\lfloor \frac{k-2pq(3d+1)}{q} \right\rfloor.
\]
By our assumption on the size of $k$ we have
\[
m := \left\lfloor \frac{k-2pq(3d+1)}{q} \right\rfloor > \lfloor (2d+2)p \rfloor.
\]
Apply Proposition \ref{fdsplit} to $\phi^{'}$ with $m$ as above.  This yields a map $\phi^{''}$
homotopic to $\phi^{'}$ which has a direct summand of the form $\oplus_{i=1}^{m-2d} e_1$.
Straightforward calculation shows that
\[
m-2d \geq \left\lfloor \frac{k-pq(8d+3)}{q} \right\rfloor > \frac{k-pq(8d+4)}{q} \geq p
\left\lfloor \frac{k-pq(8d+4)}{pq} \right\rfloor =pM.
\]
The summand $\oplus_{i=1}^{pM} e_1$ is clearly homotopic to $\oplus_{i=1}^{M} ev_{1/2}$, and the
corollary follows.
\end{proof}

\section{An extension result}

We need the following corollary of \cite[Cor. 4.6]{key-2}
\begin{cors}\label{thm:aueq}
Let $B$ be a separable, nuclear, unital, residually finite-dimensional
C*-algebra.  Let $(\pi_k)_{k=1}^\infty$
be a sequence of unital finite dimensional representations of $B$ which separates the
points of $B$ and such that each representation occurs infinitely many times.
For any unital C*-algebra $A$ and any two unital $*$-homomorphisms $\alpha,\beta:B\to M_k(A)$ with $KK(\alpha)=KK(\beta)$,
 there is a sequence of unitaries $u_n\in M_{k+r(n)}(A)$, where $r(n)$ is the
rank of the projection $\left(\pi_1\oplus\cdots\oplus \pi_n\right)(1_B)$,
such that for all $b\in B$
\[\lim_{n\to \infty} \|u_n\left(\alpha(b)\oplus\pi_1(b)\oplus\cdots\oplus
\pi_n(b)\right)u_n^*-\beta(b)\oplus\pi_1(b)\oplus\cdots\oplus \pi_n(b)\|=0.\]
 \end{cors}
For a finite dimensional representation $\pi:B \to M_k$, the induced
$*$-homomorphism $\pi\otimes 1_A:B \to M_k(A)$, $b\mapsto \pi(b)\otimes 1_A$
was  also denoted by $\pi$ in the statement above.

\begin{props}\label{UCTprop}
Let $p$ and $q$ be relatively prime integers strictly greater than one. There exists
a constant $K \in \mathbb{N}$ such that the following holds:  for any unital  C$^{*}$-algebra $A$
and unital $*$-homomorphisms $\psi, \gamma: \mathrm{I}_{p,q} \to
A$, the $*$-homomorphisms $\psi \oplus K\, ev_{x_{0}},\gamma
\oplus K \, ev_{x_{0}}:\mathrm{I}_{p,q}\to\mathrm{M}_{Kpq+1}(A)$ are homotopic for any point $x_0\in (0,1)$.
\end{props}

\begin{proof}
Let $a,b,m>0$ be  integers such that $ap+bq=mpq+1$. Throughout the proof
 $\mathrm{I}_{p,q}$ will be denoted by $B$.
 Define $*$-homomorphisms
$\alpha, \beta: B \to M_{mpq+1}(B)$ by
$\alpha=\mathrm{id}_B\oplus m\,(ev_{x_{0}}\otimes 1_{B} )$ and $\beta=a\,(e_0\otimes
1_{B})\oplus b\,(e_1\otimes 1_{B})$. Since $B$ is semiprojective, there
exist a finite subset $\mathcal{F}\subset B$ and $\delta>0$ such
that if $\mu,\nu:B\to A$ are two unital $*$-homomorphisms such that
$\|\mu(f)-\nu(f)\|<\delta$ for all $f\in \mathcal{F}$, then $\mu$
is homotopic to $\nu$. By Corollary \ref{thm:aueq} there are points
$t_1,...,t_r\in (0,1)$ such that if $\eta=(ev_{t_{1}}\otimes 1_{B})\oplus
(ev_{t_{2}}\otimes 1_{B})\oplus\cdots\oplus (ev_{t_{r}}\otimes 1_{B})$, then there is a unitary $u\in
M_{(m+r)pq+1}(B)$ such that
\[\|u\left(\alpha(f)\oplus \eta(f)\right)u^*-\beta(f)\oplus \eta(f)\|<\delta\]
for all $f\in \mathcal{F}$. By our choice of $\mathcal{F}$ and $\delta$, it
follows that $u(\alpha\oplus \eta) u^*$ is homotopic to $\beta\oplus \eta$.
Since the unitary group of  $M_{(m+r)pq+1}(B)$ is path connected \cite{key-4}
and since $\eta$ is homotopic to $r(ev_{x_0}\otimes  1_{B})$ we deduce that
$\alpha\oplus r(ev_{x_0}\otimes  1_{B})$ and
$\beta\oplus r(ev_{x_0}\otimes  1_{B})$ are homotopic as $*$-homomorphisms
from $\mathrm{I}_{p,q}$ to $M_{(m+r)pq+1}(B)$.
Consequently, for any unital $*$-homomorphism $\psi:B\to A$,
\[\left( \mathrm{id}_{(m+r)pq+1}\otimes\psi\right)\circ \left(\alpha\oplus r(ev_{x_0}\otimes
1_{B})\right)=\psi\oplus (m+r)\,(ev_{x_0}\otimes 1_A)\] is homotopic to
\[\left( \mathrm{id}_{(m+r)pq+1}\otimes\psi\right)\circ \left(\beta\oplus r(ev_{x_0}\otimes
1_{B})\right)=a(e_0\otimes 1_A) \oplus b(e_1 \otimes 1_A)\oplus r(ev_{x_0}\otimes 1_A).\]
It follows that if $\gamma:B\to A$ is any other unital $*$-homomorphism,
and $K=m+r$, then $\psi\oplus K\,ev_{x_0}$ is homotopic to $\gamma\oplus K\,ev_{x_0}$ since they are both homotopic to
$a(e_0\otimes 1_A) \oplus b(e_1 \otimes 1_A)\oplus r(ev_{x_0}\otimes 1_A)$.
\end{proof}

\begin{thms}\label{homotopic}
There is a constant $L>0$ such that the following statement holds:
Let $X$ be a compact metric space of covering dimension $d <
\infty$, and let \mbox{$k \geq L(pq)^2(d+1)$} be a natural number.
It follows that any two unital $*$-homomorphisms
\[
\phi, \psi:\mathrm{I}_{p,q} \to \mathrm{C}(X) \otimes \mathrm{M}_k
\]
are homotopic.
\end{thms}

\begin{proof}
Set $L = 16(K+1)$, where $K$ is the constant of Proposition \ref{UCTprop}.
We will first prove the theorem under the assumption that $X$ is a finite CW-complex, and
requiring only that
\[
k \geq \frac{L}{2}(pq)^2(d+1) = 8(K+1)(pq)^2(d+1).
\]
We will then use the semiprojectivity of $\mathrm{I}_{p,q}$ to deduce the general case.

Assume that $X$ is a finite CW-complex.  Using
Proposition\ref{equivalence} and the fact that $\mathrm{dim}(X) =
d$, we may assume that $N_g^{\phi}, \ N_g^{\psi} < d+1$.  Since
$\phi$ and $\psi$ now satisfy the hypotheses of Corollary
\ref{fdsplit2}, we may simply assume that both $\phi$ and $\psi$
have a direct summand of the form
\[
 \bigoplus_{i=1}^{M} ev_{1/2},
\]
where
\[
M = \left\lfloor \frac{k-pq(8d+4)}{pq} \right\rfloor.
\]
Straightforward calculation then shows that
\[
k -Mpq \leq k-pq \left\lfloor \frac{k-pq(8d+4)}{pq} \right\rfloor \leq k- [k-pq(8d+5)] \leq
8pq(d+1).
\]

Let us write
\[
\phi = \phi^{'} \oplus \left(\bigoplus_{i=1}^{M} ev_{1/2} \right); \ \
\psi = \psi^{'} \oplus \left(\bigoplus_{i=1}^{M} ev_{1/2} \right).
\]
Set $l = \mathrm{rank}(\phi^{'}(1))=\mathrm{rank}(\psi^{'}(1)) \leq 8pq(d+1)$, so that
$\phi^{'},\psi^{'}:\mathrm{I}_{p,q} \to \mathrm{C}(X) \otimes \mathrm{M}_l$ are unital
$*$-homomorphisms.  Set
\[
\phi^{''}=\phi^{'} \oplus K(1_{\mathrm{C}(X) \otimes \mathrm{M}_l} \otimes ev_{1/2}); \ \
\psi^{''}=\psi^{'} \oplus K(1_{\mathrm{C}(X) \otimes \mathrm{M}_l} \otimes ev_{1/2}).
\]
It follows from Proposition \ref{UCTprop} that
\[
\phi^{''}, \ \psi^{''}: \mathrm{I}_{p,q} \to \mathrm{C}(X) \otimes \mathrm{M}_{(Kpq+1)l}
\]
are homotopic.  Straightforward calculation shows that
\[
Kl \leq 8Kpq(d+1) \leq \left\lfloor \frac{L(pq)^2(d+1)-pq(8d+4)}{pq} \right\rfloor \leq
\left\lfloor \frac{k-pq(8d+4)}{pq} \right\rfloor =M.
\]
We may therefore view $\phi^{''}$ and $\psi^{''}$ as direct summands of $\phi$ and $\psi$,
respectively:
\[
\phi = \phi^{''} \oplus \left(\bigoplus_{i=1}^{M-Kl} ev_{1/2} \right); \ \ \psi = \psi^{''} \oplus
\left(\bigoplus_{i=1}^{M-Kl} ev_{1/2} \right).
\]
The homotopy between $\phi^{''}$ and $\psi^{''}$ now provides the desired homotopy between $\phi$
and $\psi$.

Now suppose that $X$ is only a metric space of finite covering dimension $d$.  We may embed $X$ into a bounded
subset of $\mathbb{R}^{2d+1}$, and so write $X = \cap_n X_n$, where $(X_n)$ is a decreasing sequence of
polyhedra.  By the semiprojectivity of $\mathrm{I}_{p,q}$ we have
\[
[\mathrm{I}_{p,q}, \mathrm{M}_k(\mathrm{C}(X))] = \lim_{n \to \infty} [\mathrm{I}_{p,q}, \mathrm{M}_k(\mathrm{C}(X_n))].
\]
We may therefore assume that the homotopy classes of our given maps $\phi$ and $\psi$ lie in some
$[\mathrm{I}_{p,q}, \mathrm{M}_k(\mathrm{C}(X_n))]$.  Having proved the theorem for finite CW-complexes,
we conclude that $\phi$ and $\psi$ are homotopic if
\[
k \geq \frac{L}{2}(pq)^2(\mathrm{dim}(X_n)+1) \geq \frac{L}{2}(pq)^2(2d+2) = L(pq)^2(d+1).
\]
This proves the theorem proper.
\end{proof}

\begin{cors}\label{extend2}
Let $p$ and $q$ be relatively prime positive integers strictly greater than one,
and let $X$, $L$, and $k$ be as in Theorem \ref{homotopic}.  Suppose that $Y \subseteq X$
is closed, and that we are given a unital $*$-homomorphism
$\phi:\mathrm{I}_{p,q} \to \mathrm{C}(Y) \otimes \mathrm{M}_k.$
It follows that there is a unital $*$-homomorphism $\psi:\mathrm{I}_{p,q} \to
\mathrm{C}(X) \otimes \mathrm{M}_k$ such that $\psi|_Y = \phi$.
\end{cors}

\begin{proof}
By the semiprojectivity of $\mathrm{I}_{p,q}$ we can extend $\phi$
to the closure of some open neighbourhood $O$ of $Y$, i.e., we may
assume that $\phi:\mathrm{I}_{p,q} \to \mathrm{C}(\overline{O})
\otimes \mathrm{M}_k$ without changing the original definition of
$\phi$ over $Y$. As explained in  subsection ~\ref{2.2}, $F^k$ is
nonempty. Choose a point $\gamma \in F^k$.
 By Theorem \ref{homotopic}, $\phi$ and $1_{\mathrm{C}(\overline{O})} \otimes
\gamma$ are homotopic as maps from $\mathrm{I}_{p,q}$ into $\mathrm{C}(\overline{O})
\otimes \mathrm{M}_k$.  Let us denote this homotopy by
\[H:[0,1] \to \mathrm{Hom}_1(\mathrm{I}_{p,q};
\mathrm{C}(\overline{O}) \otimes \mathrm{M}_k),
\]
where $H(0)=\phi$.

Find a continuous map $f:X \to [0,1]$ which is equal to zero on $Y$
and equal to one off $O$.  Define $\psi:\mathrm{I}_{p,q} \to \mathrm{C}(X) \otimes
\mathrm{M}_k$ by the formula
\[
\psi_x = \left\{ \begin{array}{ll} H(f(x))_x, & x \in \overline{O} \\
\gamma, x \notin \overline{O} \end{array} \right.
\]
One checks that $\psi$ so defined has the required property.
\end{proof}

\section{Maps from $\mathrm{I}_{p,q}$ into recursive subhomogeneous algebras}

\subsection{Recursive subhomogeneous algebras}\label{rsh}
Let us recall some of the terminology and results from \cite{ph1}.
\begin{dfs}\label{rshdef}

\begin{enumerate}
    \item[(i)] if $X$ is a compact Hausdorff space and $k \in \mathbb{N}$, then $\mathrm{M}_k(\mathrm{C}(X))$
    is a recursive subhomogeneous algebra (RSH algebra);
    \item[(ii)] if $A$ is a recursive subhomogeneous algebra, $X$ is a
    compact Hausdorff space, $X^{(0)} \subseteq X$ is closed, $\phi:A \to \mathrm{M}_k(\mathrm{C}(X^{(0)}))$
    is a unital $*$-homomorphism, and
    $\rho:\mathrm{M}_k(\mathrm{C}(X)) \to \mathrm{M}_k(\mathrm{C}(X^{(0)}))$ is the restriction homomorphism, then
    the pullback
    \[
    A \oplus_{\mathrm{M}_k(\mathrm{C}(X^{(0)}))} \mathrm{M}_k(\mathrm{C}(X))
    =\{(a,f) \in A \oplus \mathrm{M}_k(\mathrm{C}(X)) \ | \
    \phi(a) = \rho(f) \}
    \]
    is a recursive subhomogeneous algebra.
\end{enumerate}
\end{dfs}
\noindent It is clear from the definition above that a C$^*$-algebra $R$ is an RSH algebra if and
only if it can be written in the form
\begin{equation}\label{decomp2}
R = \left[ \cdots \left[  \left[ C_0 \oplus_{C_1^{(0)}} C_1 \right] \oplus_{C_2^{(0)}} C_2 \right]
\cdots \right] \oplus_{C_l^{(0)}} C_l,
\end{equation}
with $C_k = \mathrm{M}_{n(k)}(\mathrm{C}(X_k))$ for compact Hausdorff spaces $X_k$ and integers
$n(k)$, with $C_k^{(0)}=\mathrm{M}_{n(k)}(\mathrm{C}(X_k^{(0)}))$ for compact subsets $X_k^{(0)}
\subseteq X$, and where the maps $C_k \to C_k^{(0)}$ are always the restriction maps.  We call the
C$^*$-algebra
\[
R_k = \left[ \cdots \left[  \left[ C_0 \oplus_{C_1^{(0)}} C_1 \right] \oplus_{C_2^{(0)}} C_2
\right] \cdots \right] \oplus_{C_k^{(0)}} C_k
\]
the {\it $k^{\mathrm{th}}$ stage algebra of $R$}.  Let $\mathrm{Prim}_n(R)$ denote the space of
irreducible representations of $R$ of dimension $n$. We say that an RSH algebra has finite
topological dimension if $\mathrm{dim}(\mathrm{Prim}_n(R))$ is finite for each $n$;  if $R$ has
finite topological dimension, then we call $d:=\mathrm{max}_n \mathrm{dim}(\mathrm{Prim}_n(R))$
the {\it topological dimension of $R$}.  If $R$ is separable, then the $X_k$ can be taken to be
metrisable (\cite[Proposition 2.13]{ph1}).  Finally, if $R$ has no irreducible representations of
dimension less than or equal to $N$, then we may assume that $n(k)
> N$.  We refer to the smallest of the $n(k)$ as the {\it minimum matrix size of $R$}.

\subsection{An existence theorem}

\begin{thms}\label{ipqmap}
Let $R$ be a separable RSH algebra of finite topological dimension
$d$ and minimum matrix size $n$.  Let $p$ and $q$ be relatively
prime integers strictly greater than one.  Suppose that $n \geq
L(pq)^2(d+1)$, where $L$ is the constant of Theorem
\ref{homotopic}.  It follows that there is a unital
$*$-homomorphism $\gamma:\mathrm{I}_{p,q} \to R$.
\end{thms}

\begin{proof}
We proceed by induction on the index $k$ from subsection \ref{rsh}. We have a decomposition
\[
R = \left[ \cdots \left[  \left[ C_0 \oplus_{C_1^{(0)}} C_1 \right] \oplus_{C_2^{(0)}} C_2 \right]
\cdots \right] \oplus_{C_l^{(0)}} C_l
\]
as in (\ref{decomp2}) above, where $C_0 =
\mathrm{M}_{n(0)}(\mathrm{C}(X_0))$ and $n(0) \geq n > pq$. As
explained in subsection~\ref{2.2}, $F^{n(0)}$ is not empty. Choose
$\psi \in F^{n(0)}$. It follows that $\gamma_0 :=
1_{\mathrm{C}(X_0)} \otimes \psi$ defines a unital
$*$-homomorphism from $\mathrm{I}_{p,q}$ into $C_0$.

Suppose $k<l$, and that we have found a unital $*$-homomorphism $\gamma_k:\mathrm{I}_{p,q} \to
R_k$.  We will prove that $\gamma_k$ can be extended to a unital $*$-homomorphism
$\gamma_{k+1}:\mathrm{I}_{p,q} \to R_{k+1}$.  Starting with $\gamma_0$ and applying this inductive
result repeatedly will yield the map $\gamma$ required by the theorem.  We have
\[
R_{k+1} = R_k \oplus_{C_{k+1}^{(0)}} C_{k+1}.
\]
Notice that $\gamma_k$ defines a unital $*$-homomorphism from
$\mathrm{I}_{p,q}$ into $R_k \oplus C_{k+1}^{(0)}$ in a natural
way---the map into the summand $R_k$ is simply $\gamma_k$ itself,
while the map into the summand $C_{k+1}^{(0)}$ is the composition
of $\gamma_k$ with the clutching map $\phi:R_k \to C_{k+1}^{(0)}$
(cf. Definition \ref{rshdef}).  Viewing $R_{k+1}$ as a subalgebra
of $R_k \oplus C_{k+1}$, we see that our task is simply to extend
the map $\phi \circ \gamma_k:\mathrm{I}_{p,q} \to C_{k+1}^{(0)}$
to all of
$C_{k+1}=\mathrm{M}_{n(k+1)}(\mathrm{C}(X_{k+1}^{(0)}))$.  Since
$n(k+1) \geq n \geq L(pq)^2(d+1)$, the existence of the desired
extension follows from Corollary \ref{extend2}.
\end{proof}

\subsection{Proof of Theorem \ref{main}}

\begin{proof}
Let $A$ be a unital separable C$^*$-algebra.
By \cite[Proposition 6.3]{t1}, it will suffice to prove that for any $\mathrm{I}_{p,q}$, there is
a unital $*$-homomorphism $\gamma:\mathrm{I}_{p,q} \to A^{\otimes \infty}$.

By hypothesis, there is a unital subalgebra $S$ of $A^{\otimes
\infty}$ which is separable, subhomogeneous, and has no
characters.  By the main result of \cite{nw}, S is the limit of an
inductive system $(R_i,\phi_i)$, where each $R_i$ is a (unital)
separable RSH algebra of finite topological dimension and each
$\phi_i$ is injective and unital. Suppose, contrary to our desire,
that $X_i := \mathrm{Prim}_1(R_i) \neq \emptyset$ for each $i \in
\mathbb{N}$.. Each $\phi_i:R_i \to R_{i+1}$ induces a continuous
map $\phi_i^{\sharp}:X_{i+1} \to X_i$.  Since each $X_i$ is
compact, the limit of the inverse system
$(X_i,\phi_{i-1}^{\sharp})$ is not empty.  In other words, there
is an element of $X_1$ which has a pre-image in $X_{i+1}$ under
each composed map $\phi_1^{\sharp} \circ \cdots \circ
\phi_{i}^{\sharp}$.  It follows that $S$ has a character, contrary
to our assumption.  We therefore conclude that $R_i$ has no
characters for some $i \in \mathbb{N}$.  Since the $\phi_i$ are
injective, we have that $A^{\otimes\infty}$ contains, unitally, a
recursive subhomogeneous algebra $R:=\phi_{i \infty}(R_i)$ of
finite topological dimension which has no characters.

Let $n >1$ be the minimum
matrix size of $R$.  Find a natural number $m$ such that $n^m/(md+1) > L(pq)^2$.  It follows from
\cite[Proposition 3.4]{ph1} that the topological dimension of $R^{\otimes m}$ is at most $md$,
while the minimum matrix size of $R^{\otimes m}$ is at least $n^m$.  Applying Theorem \ref{ipqmap}
we obtain a unital $*$-homomorphism
\[
\gamma:\mathrm{I}_{p,q} \to R^{\otimes m} \hookrightarrow (A^{\otimes \infty})^m \cong A^{\otimes
\infty},
\]
as required.
\end{proof}

\subsection{Examples} Let us explain now why the examples (a)-(f) of the introduction satisfy the hypotheses of
Theorem \ref{main}

\begin{enumerate}
\item[(a)]  Let $A$ be a unital simple separable exact C$^{*}$-algebra containing an infinite projection.  It follows from
a result of Kirchberg (\cite{k}) that $A^{\otimes \infty}$---even $A^{\otimes 2}$---is purely infinite and simple,
and so has real rank zero.  By Proposition 5.7 of \cite{pr}, there is a unital $*$-homomorphism $\phi:F \to
A^{\otimes \infty}$,
where $F$ is a finite-dimensional (hence subhomogeneous) C$^*$-algebra without characters.
\item[(b)] Let $A$ be a unital separable ASH algebra without characters.  Following the arguments in the proof of
Theorem \ref{main}, we see that there must be a unital subalgebra of $A$ which is subhomogeneous without characters.
\item[(c)] If $A$ is properly infinite, then there is a unital
embedding of $\mathcal{O}_{\infty}$ into $A$;  $\mathcal{O}_{\infty}$ is $\mathcal{Z}$-stable
by the Kirchberg-Phillips classification, and so any $\mathrm{I}_{p,q}$
embeds unitally into $A$.
\item[(d)] Let $A$ be a unital separable C$^*$-algebra of real rank zero.  By Proposition 5.7 of \cite{pr}, there is a unital map
$\phi:F \to A$, where $F$ is a finite-dimensional C$^*$-algebra without characters.
\item[(e)] Let $X$ be a compact infinite metric space, and $\alpha:X \to X$ a minimal homeomorphism.  It follows from Theorem
2.7 of \cite{lp} that the crossed product $\mathrm{C}^*(X,\mathbb{Z},\alpha)$ contains a recursive subhomogeneous algebra
without characters.
\item[(f)] There are several examples which show that Banach algebra $\mathrm{K}$-theory and traces do
not form a complete invariant for simple unital separable amenable C$^*$-algebras.  The first of these is
due to R\o rdam, and consists of a simple unital separable amenable C$^*$-algebra $A$ containing both a finite and an
infinite projection.  By the theorem of Kirchberg cited in (a) we have that $A \otimes A$ is purely infinite, and
so following the arguments of (a) we see that $A$ satisfies the hypotheses of Theorem \ref{main}.  Other
examples were produced by the second name author in \cite{t2} and \cite{t3}.  These algebras are ASH and non-type-I,
and so satisfy the hypotheses of Theorem \ref{main} by the arguments of (b) above.
\end{enumerate}

\section{Concluding remarks}

One consequence of Theorem \ref{ipqmap} is the following:

\begin{cors}\label{ipqash}
Let $A$ be a simple separable unital non-type-I inductive limit of RSH algebras with
slow dimension growth.  Given any relatively prime integers $p,q>1$, there is a unital
$*$-homomorphism $\phi:\mathrm{I}_{p,q} \to A$.
\end{cors}

\noindent
This result may be viewed as a step toward proving that an algebra $A$ as in Corollary \ref{ipqash}
is $\mathcal{Z}$-stable.  To get the stronger conclusion one needs to parlay the map $\phi$ into
an approximately central sequence of $*$-homomorphisms $\phi_n:\mathrm{I}_{p,q} \to A$.
This improvement, however difficult it may be to realise technically, is at least eminently reasonable
for reasons of spectral multiplicity.  Given an inductive sequence $(A_i,\gamma_i)$ of RSH algebras
with slow dimension growth (see \cite{ph1} for the definition of this property), one knows that the
commutant of the image of $A_i$ in $A_j$ over each point in the spectrum of $A_j$ is well-approximated
by a finite-dimensional C$^*$-algebra $F$, all of whose simple summands have large dimension compared with
the dimension of the spectrum of $A_j$.  It stands to reason that one should be able to map  an RSH algebra, all of whose matrix
fibres are large in comparison with its topological dimension, into the approximate
commutant of the image of a finite set $F \subseteq A_i$ in $A_j$.   An application of Theorem \ref{ipqmap} will
then provide a $*$-homomorphism $\phi:\mathrm{I}_{p,q} \to A_j$ which almost commutes with the image of
$F$ in $A_j$, leading to the $\mathcal{Z}$-stability of $A = \lim_{i \to \infty}(A_i,\gamma_i)$.

Proving that $A$ as in Corollary \ref{ipqash} is $\mathcal{Z}$-stable
would be an extremely important step toward the confirmation
of Elliott's classification conjecture for simple unital separable ASH algebras of slow dimension
growth, and would, given Winter's recent classification theorem (\cite{w1}) and the decomposition results of
Lin and Phillips (\cite{lp}), already confirm Elliott's conjecture for the class of all
crossed product C$^*$-algebras arising from the action of a minimal uniquely ergodic diffeomorphism
with smooth inverse on a compact manifold.

\end{document}